\documentclass[12pt,a4paper,draft]{article}
\usepackage[russian]{babel}
\usepackage{tocvsec2}
\usepackage{amsmath}
\usepackage{amsfonts}
\usepackage{amssymb}
\usepackage{graphicx}
\usepackage{pgf,tikz}
\usetikzlibrary{arrows}
\usepackage[labelformat=simple]{subcaption}

\newcounter{example}
\newcounter{tesis}
\newcounter{utv}
\newenvironment{opr}{\smallskip\refstepcounter{example}\textbf{Определение \arabic{example}. }}{\par}

\newenvironment{teor}{\smallskip\par\refstepcounter{tesis}\textbf{Теорема \arabic{tesis}. }}{\par}

\newenvironment{lemma}{\smallskip\refstepcounter{utv}\textbf{Лемма \arabic{utv}.}}{\par}

\newenvironment{dokvo}{\smallskip\par\textit{Доказательство. }}{$\square$\par}

\begin{document}
\Large
\thispagestyle{empty}

%\vspace{-1cm}
%\begin{center}
%НАЦИОНАЛЬНЫЙ ИССЛЕДОВАТЕЛЬСКИЙ УНИВЕРСИТЕТ\\
%ВЫСШАЯ ШКОЛА ЭКОНОМИКИ\\[1ex]
%Факультет математики\\[9ex]

%Курсовая работа\\[8ex]
%{\bf ПОКРЫТИЕ ПРЯМОУГОЛЬНИКА КВАДРАТАМИ}\\[9ex]
%\end{center}
%\hangindent=7,5cm
%\hangafter=0
%\noindent
%Выполнил:\\
%бакалавр 1-го года обучения\\
%{\it Ожегов Федор Юрьевич}\\[6ex]
%Научный руководитель:\\
%кандидат физ.-мат. наук, доцент\\
%{\it Скопенков Михаил Борисович}\\[8ex]
%\begin{center}
%Москва\\
%2018
%\end{center}
\small
\begin{center}
 Оклеивание прямоугольника квадратами \par
Ожегов Федор Юрьевич \par 
 Национальный исследовательский университет Высшая школа экономики
\end{center}

\footnotesize
Аннотация. В работе доказывается, что оклеить прямоугольник $1 \times b$  равными квадратами с двух сторон в один слой можно тогда и только тогда, когда $b$ = $p \pm \sqrt{p^2 - r^2}$, где $p \geq r \geq 0$ рациональны.
\bigskip
\bigskip

\small

 В данной работе решена следущая задача: \par
Дан конверт в форме прямоугольника $a \times b$. При каких вещественных $a$ и $b$ его можно оклеить квадратными марками одинакового размера без просветов и наложений с обеих сторон? Квадраты разрешается перегибать через край прямоугольника. \par

\begin{teor} 
 Оклеить прямоугольник $1 \times b$  равными квадратами с двух сторон в один слой можно тогда и только тогда, когда $b$ = $p \pm \sqrt{p^2 - r^2}$, где $p \geq r \geq 0$ рациональны.
\end{teor}

Понятие оклеивания формально определяется в определении 1 ниже.

Рассматриваемая задача похожа на задачу о покрытии квадратами двумерных поверхностей из работы \cite{kenyon}.

В ходе доказательства теоремы 1 мы используем идею из работы \cite{covering} о сопоставлении оклеиванию конверта периодического замощения плоскости и идею из работы \cite{kvant} о представлении этого замощения в виде объединения счетного числа полос квадратов. 

Начнем с неформальной мотивировки.
Покрасим стороны конверта в два цвета, положим его на плоскость и начнем его перекатывать через стороны. Таким образом получим прямоугольную решетку $\Lambda$ на плоскости. Оклеивание конверта даст периодическое замощение плоскости, переходящее в себя при центральных симметриях относительно узлов решетки $\Lambda$. Тогда исходная задача сводится к поиску таких замощений.

\begin{opr}
Зафиксируем декартову систему координат  на плоскости.
 Решетку $\Lambda  =  \{(mb , n): m,n \in \mathbb{Z} \}$  на плоскости назовем \emph{порожденной} прямоугольником $1 \times b$.
Введем отношение эквивалентности: назовем две точки плоскости эквивалентными, если одну можно перевести в другую композицией центральных симметрий относительно узлов решетки $\Lambda$.
\par

\emph{Оклеиванием} (с двух сторон в один слой) прямоугольника $1 \times b$ равными квадратами назовем такое конечное множество $G$ непересекающихся по внутренним точкам равных квадратов на плоскости, что:
\begin{itemize}
\item любая точка плоскости эквивалентна некоторой точке объединения квадратов множества $G$;
\item  никакие две точки внутренности объединения квадратов из множества $G$ не эквивалентны.
\end{itemize}
\end{opr}

Такие определения помогают сразу рассматривать оклеивание конверта квадратами, как некий объект на плоскости.

\begin{opr}\label{opr5}
\emph{Полосой квадратов} назовем такую бесконечную в обе стороны последовательность квадратов на плоскости, что любые два соседних члена этой последовательности имеют общую сторону и любые два несоседних члена последовательности не имеют общих вершин. 
\end{opr}

\addtocontents{toc}{\vspace{2mm}}

\begin{lemma}\label{lem1}
Рассмотрим оклеивание $G$, состоящее из квадратов  $K_1$ \dots $K_g$. Рассмотрим множество $M_j$ квадратов на плоскости, получаемых из квадрата $K_j$ композициями центральных симметрий относительно узлов решетки $\Lambda$. Тогда $\bigcup^g_{i = 1}\bigcup_{K \in M_i}K$ является замощением плоскости квадратами. При этом оно переходит в себя при любой композиции центральных симметрий относительно узлов решетки $\Lambda$.

\begin{dokvo}
В множестве $G$ присутствуют представители всех классов эквивалентности. Тогда для любой точки плоскости существует эквивалентная ей точка оклеивания, а значит, $\mathbb{R}^2 = \bigcup^g_{i = 1}\bigcup_{K \in M_i}K$.
Покажем, что любые два квадрата из $\bigcup^g_{i = 1} M_i$ либо совпадают, либо не имеют общих внутренних точек.
Пусть квадраты $K$ и $K'$ из $\bigcup^g_{i = 1} M_{i}$ имеют общую внутреннюю точку. Так как внутренние точки несовпадающих множеств $M_i$ и $M_j$ являются представителями разных классов эквивалентности, то  $K$ и $K'$ лежат в одном множестве $M_l$. Тогда, если $K$ не совпадает с $K'$, то в $K$ найдутся две внутренних эквивалентных точки, что невозможно из определения оклеивания. Тогда  $\mathbb{R}^2$ представлено в виде объединения счетного числа непересекающихся по внутренним точкам квадратов.
 \par
 Так как при композиции центральных симметрий относительно узлов решетки $\Lambda$ квадраты, принадлежащие $M_i$, переходят в квадраты, принадлежащие
$M_i$, то и замощение переходит в себя при композиции центральных симметрий относительно узлов решетки $\Lambda$.
\end{dokvo}
\end{lemma}

\begin{lemma}\label{lem2}
 Замощение из леммы 1 является объединением счетного числа полос квадратов.

\textit{Доказательство.} 
Будем обозначать вершины квадратов, начиная с левого нижнего угла, идя по часовой стрелке. Рассмотрим два случая: \par
 Случай $(1)$ : найдутся два квадрата $A_0B_0C_0D_0$ и $A_1B_1C_1D_1$ замощения из леммы 1, пересекающиеся только по части стороны (но не только по вершине). 
Без ограничения общности можно считать, что $A_1 \in B_0C_0, C_0 \in A_1D_1$. 

Ясно, что точки $A_0 , B_1$ лежат по разные стороны от прямой $A_1C_0$. Тогда,  так как  все квадраты замощения равны, то квадрат $D_0C_0C_2D_2$, прилегающий к углу $D_1C_0D_0$ восстановится однозначно, так как одной стороной он прилегает к прямой $C_0D_0$, а другой - к $A_1D_1$. Аналогично построим квадрат, прилегающий к стороне $C_0C_2$ квадрата $D_0C_0C_2D_2$ и стороне $C_1D_1$ квадрата $A_1B_1C_1D_1$. Таким образом построим полосу $L_0$ квадратов, содержащую $A_0B_0C_0D_0$, и полосу $L_1$ квадратов, содержащую $A_1B_1C_1D_1$. Теперь рассмотрим такой квадрат $A_2B_2C_2D_2$, что $A_2 \in B_1C_1$. Ясно, что квадрат, прилегающий к стороне $C_2D_2$ квадрата $A_2B_2C_2D_2$, восстановится однозначно, а значит, восстановится и полоса $L_2$, содержащая квадрат $A_2B_2C_2D_2$. Аналогично построим счетное число полос.\par
 Случай $(2)$: любые два квадрата замощения, имеющих хотя бы две общие точки, имеют общую сторону. Ясно, что тогда равные квадраты образуют квадратную сетку на плоскости, в частности, образуют объединение счетного числа полос. $\square$

\end{lemma}

\hspace{2 mm} \textit{Доказательство теоремы 1.} \emph{Необходимость}. Рассмотрим некоторое оклеивание $G$, состоящее из $g$ равных квадратов. Тогда, из леммы 1, существует некоторое замощение плоскости квадратами, полученное из оклеивания $G$. В силу леммы 2 это замощение является объединением счетного числа полос квадратов. 
 Пусть $\alpha$ - такой угол, что вектор $\overline{n} = (-\sin\alpha, \cos\alpha )$ перпендикулярен направлению полос. Так как замощение переходит в себя при переносах на векторы вида $(2bt_1, 2t_2)$, где $t_1, t_2 \in \mathbb{Z}$, то проекции векторов $(-2b, 0)$ и $(0, 2)$ на $\overline{n}$ имеют длину, кратную стороне $a$ квадрата. 
Обозначим $2b\sin\alpha  =: aq_1, 2\cos\alpha =: aq_2$, где $q_1, q_2 \in \mathbb{Z}$. Тогда из равенства $(\sin\alpha)^2 + (\cos\alpha)^2 = 1$ получаем $a^2(\frac{q^2_1}{4b^2} + \frac{q^2_2}{4}) = 1$. Заметим, что прямоугольник $1 \times 2b$ содержит представителей всех классов эквивалентности и не содержит внутри себя эквивалентных точек. Тогда площадь прямоугольника $1 \times 2b$ равна площади оклеивания, то есть $ga^2$. Если $q_2 = 0$, то $2b = ga^2 = g(\frac{2b}{q_1})^2$, а значит, $b$ рационально. Тогда можно взять $r = 0, p = \frac{b}{2} $. Если же $q_2 \neq 0$, то, так как $a^2 = \frac{b}{g}$, получаем $b^2 - b\frac{4g}{q^2_2} + (\frac{q_1}{q_2})^2 = 0$. Обозначив $p:= \frac{2g}{q^2_2}$, $r:= (\frac{q_1}{q_2})^2$ получаем $b = p \pm \sqrt{p^2 - r^2}$, где $r,p \in \mathbb{Q}$.

\hspace{2 mm} \emph{Достаточность}. Пусть известно, что $b = p \pm \sqrt{p^2 - r^2}$, где $p \geq r \geq 0$ рациональны. Пусть $r = \frac{n}{m}$ , где $m,n \geq 0$ - взаимнопростые целые числа. Рассмотрим прямоугольник $K$ с вершинами в точках $(bm, n)$, $(bm + \frac{bn}{b^2m^2 + n^2}, n - \frac{b^2m}{b^2m^2 + n^2})$ , $(-bm + \frac{bn}{b^2m^2 + n^2}, -n - \frac{b^2m}{b^2m^2 + n^2})$, $(-bm, -n)$. Стороны прямоугольника $K$ равны $2\sqrt{b^2m^2 + n^2} $ и $\frac{b}{\sqrt{b^2m^2 +n^2}}$, а их отношение - $2\frac{b^2m^2 + n^2}{b} = 2\frac{b^2 + r^2}{b}m^2 = 2\frac{2p^2 \pm 2\sqrt{p^2 - r^2}}{p \pm \sqrt{p^2 - r^2}}m^2 = 4m^2p \in \mathbb{Q}$. Разрежем $K$ на некоторое количество квадратов и получим множество $G$ непересекающихся по внутренним точкам квадратов. Докажем, что $G$ - искомое оклеивание. Обозначим через $K'$ объединение прямоугольника $K$ и его симметрии относительно точки $(0,0)$. \par
Покажем, что никакие две внутренние точки прямоугольника $K$ не эквивалентны (а значит, никакие две внутренние точки объединения квадратов из множества $G$ тоже не эквивалентны). Предположим обратное: существуют две внутренние эквивалентные точки прямоугольника $K$. Тогда они либо совмещаются параллельным переносом  на вектор вида $(2bt_1,2t_2)$, где $t_1, t_2 \in \mathbb{Z}$ не равны нулю одновременно, либо центрально симметричны относительно некоторого узла $A$ решетки $\Lambda$. В последнем случае рассмотим композицию центральной симметрии относительно узла $A$ с симметрией относительно точки $(0,0)$. Эта композиция - снова параллельный перенос на вектор вида $(2bt_1,2t_2)$, который совмещает две внутренних точки уже прямоугольника $K'$. Рассмотрим два случая. \par

Случай $(1)$:
пусть $mt_2 \neq nt_1$. Проекция вектора $(2bt_1, 2t_2)$ на сторону прямоугольника $K'$, перпендикулярную прямой $y = \frac{xn}{bm}$, равна
 $$\left| 2\left(bt_1, t_2\right) \cdot \left(-\frac{n}{\sqrt{b^2m^2+n^2}}, \frac{bm}{\sqrt{b^2m^2 + n^2}}\right) \right|= 2b\frac{\left| mt_2 - nt_1 \right|}{\sqrt{b^2m^2 + n^2}} \geq \frac{2b}{\sqrt{b^2m^2 + n^2}},$$ что не меньше длины этой стороны при $mt_2 \neq nt_1$.	
 \par
Случай $(2)$:
Пусть $mt_2 = nt_1$. Так как $m$ и $n$ взаимнопросты, то $t_1=km, t_2=kn$, где $k \in \mathbb{Z}$. Тогда проекция вектора $(2bt_1, 2t_2)$ на сторону прямоугольника $K'$, параллельную прямой  $y = \frac{xn}{bm}$, равна  $$\left| 2\left(bkm, kn\right) \cdot \left(\frac{mb}{\sqrt{b^2m^2+n^2}}, \frac{n}{\sqrt{b^2m^2 + n^2}}\right) \right| = 2|k|\sqrt{bm^2 + n^2} \geq 2\sqrt{bm^2 + n^2}, $$ что не меньше этой стороны при $k \neq 0$. Тогда внутренние точки прямоугольника $K$ или $K'$ переносом на такой вектор совместить нельзя. В обоих случаях $(1)$ и $(2)$ получаем противоречие, значит никакие две внутренние точки прямоугольника $K$ не эквивалентны.

Заметим, что площадь прямоугольника $K$ равна $2b$, что равно площади прямоугольника $1 \times 2b$, в котором присутствуют все классы эквивалентности. Тогда в $K$ присутствуют представители всех классов эквивалентности. А значит, $G$ - оклеивание. $\square$

 Естественный вопрос состоит в обобщении теоремы 1 на случай неравных квадратов. В работе \cite{kenyon} рассматривались периодические покрытия плоскости квадратами. Как мы увидели в ходе доказательства нашего критерия, оклеиванию  было сопоставлено замощение плоскости, переходящее в себя при центральных симметриях относительно узлов некоторой решетки. В частности, замощение было периодическим, поэтому, по аналогии с $\cite[Theorem 10]{kenyon}$ предлагается следущая гипотеза. \par
\textbf{Гипотеза.} Оклеить конверт $1 \times b$ квадратами с двух сторон в один слой можно тогда и только тогда, когда $b = p \pm \sqrt{q^2 - r^2}$, где $p, r, q$ рациональны и  $ q \geq p \geq \sqrt{q^2 - r^2}$.

\end{document}